\documentclass[12pt,a4paper]{article}
\usepackage{amsfonts}
\usepackage{mathrsfs}
\usepackage{amsmath}
\usepackage{amssymb}
\usepackage{booktabs}

\renewcommand{\baselinestretch}{1}
\renewcommand{\thesection}

\begin{document}
\renewcommand\refname{\centerline{References}}
\title{{Solutions for Neumann boundary value problems involving $\big(p_{1}(x), p_{2}(x)\big)$-Laplace operators
\thanks{Research supported by the National Natural Science Foundation of China
(NSFC 10971088 and NSFC 10971087).}}}
\author{Duchao Liu\footnote{Email address: liuduchao@yahoo.cn}, Xiaoyan Wang\footnote{Email address: wang264@indiana.edu} and Jinghua Yao\footnote{Email address: yaoj@indiana.edu}
 \\
{\small School of Mathematics and Statistics, Lanzhou
University, Lanzhou, 730000, China} \\
{\small Department of Mathematics, Indiana University Bloomington, IN, 47405, USA} \\
\small Department of Mathematics, Indiana University Bloomington,
IN, 47405, USA}
\date{}
\maketitle

\begin{center}
\begin{minipage}{120mm}
\baselineskip 12pt {\small {\normalsize\textbf{Abstract.}} In this paper we study the nonlinear Neumann boundary value problem of the following equations
\begin{align*}
-\text{div}(|\nabla u|^{p_{1}(x)-2}\nabla u)-\text{div}(|\nabla u|^{p_{2}(x)-2}\nabla u)+|u|^{p_{1}(x)-2}u+|u|^{p_{2}(x)-2}u=\lambda f(x,u)
\end{align*}
in a bounded smooth domain $\Omega\subset\mathbb{R}^{N}$ with
Neumann boundary condition given by
\begin{align*}
|\nabla u|^{p_{1}(x)-2}\frac{\partial u}{\partial\nu}+|\nabla u|^{p_{2}(x)-2}\frac{\partial u}{\partial\nu}=\mu g(x,u)
\end{align*}
on $\partial\Omega$. Under appropriate conditions on the source and boundary nonlinearities, we obtain a number of results on existence and multiplicity of solutions by variational methods in the framework of variable exponent Lebesgue and Sobolev spaces. \\[6pt]
\textbf{Keywords} Critical points; Variable exponent Lebesgue-Sobolev Space; Mountain-Pass Lemma; Fountain Theorem; Dual Fountain Theorem.\\[6pt]
\textbf{Mathematics Subject Classification (2000)} 35B38 35D05
35J20}
\end{minipage}
\end{center}

\section*{0. Introduction}

In the theory of electrorheological fluids, image restoration (see \cite{1,2,3} 
and the references therein) and mixture of two power-law fluids (see \cite{20}), the $p(x)$-Laplace operator,
defined as $\Delta_{p(x)}u(x):=\text{div}\big(|\nabla
u(x)|^{p(x)-2}\nabla u\big)$, plays an important role. 
After simplifications, those problems are reduced to study differential equations involving the $p(x)$-Laplace operator. For problems involving different growth rates depending
on the underlying domains, they also involve equations with $(p(x), q(x))-$ growth conditions where several $p(x)-$Laplace operators involved, interacting with one another. 
This $(p(x), q(x))-$growth condition is a natural generalization of  the anisotropic (p,q) growth condition in electrorheological fluids (see \cite{21, 22, 23}). With the development of the theory
on variable exponent Lebesgue and Sobolev spaces, lots of mathematical works have appeared and there are many interesting applications in fluids and electrorheological fluids (see the recent monograph [4] and references) in the Orlicz-type spaces framework.

Driven by  the above mentioned physical models and potential further applications, the study of differential equations
involving $p(x)$-Laplace operators has been a very interesting and
exciting topic in recent years (see in particular these nice references
\cite{4,5, 20, 21, 22, 23}).  As we will see in the following parts of the paper, for problems involving
nonstandard growth conditions, the growth conditions are crucial for the solvability of the problems and closely related to the
space structures under which we seek solutions.  During the development of last several decades, it turns out that the use of
variational methods in dealing with these problems is a far-reaching field. Many results have been
obtained on this kind of problems, for example
\cite{6,7,8,9,10,11,12,13} and so on.

In this paper, continuing our former investigations on these topics
(see \cite{6,7,8,11,12,13,14,15,16}) and functional properties of $(p_{1}(x),p_{2}(x))$-Laplace operators (see [16]),  we study the
following $(p_{1}(x), p_{2}(x))$-Laplace Neumann boundary value problem
\begin{align*}
(P)\left\{ \begin{array}{rcl}
         &-\Delta_{p_1(x)}u-\Delta_{p_2(x)}u+|u|^{p_{1}(x)-2}u+|u|^{p_{2}(x)-2}u=\lambda f(x,u), &  \mbox{ in }
         \Omega; \\
         &|\nabla u|^{p_{1}(x)-2}\frac{\partial u}{\partial\nu}+|\nabla u|^{p_{2}(x)-2}\frac{\partial u}{\partial\nu}=\mu g(x,u),&\mbox{ on } \partial\Omega,
          \end{array}\right.
\end{align*}
where $\Omega\subset\mathbb{R}^{N}$ is a bounded smooth domain and
$p_{i}(x)\in C(\overline{\Omega})$ with $p_{i}(x)>1$ for any
$x\in\overline{\Omega}$ and for $i=1, 2$; $\lambda,\mu\in\mathbb{R}$ such that $\lambda^{2}+\mu^{2}\not=0$. 

By variational arguments, we show among others that the problem admits mountain pass solution, fountain solutions and dual fountain solutions with appropriate energy behavior when multiple solutions exist. 
To deal with the boundary nonlinear term, we shall use the trace embedding theorem for our working space; To deal with the nonlinear non-homogeneous $(p_1 (x), p_2(x))-$Laplace operator and to show 
the functionals involved in our problems verifies compactness conditions of $(PS)$ type, we shall use the properties that the corresponding nonlinear function of these operator induces homeomorphism between space pairs and admits Frechet derivative of type $(S+)$ proved in our former work.  In particular, different to the use of fountain theorem, we need to show $(PS^*_c)$-type compactness conditions the proof of which is subtle  in order to apply the dual fountain theorem. See section 3 for the assumptions and statements of results.

This paper is organized as follows. We  begin by recalling the definitions of the variable
exponent Lebesgue-Sobolev spaces which can be regarded as a special
class of generalized Orlicz-Sobolev spaces and introduce some basic
properties of these spaces. For the
convenience of the readers and to make preparation for the proof in Section 3, we give some properties
of the $(p_{1}(x), p_{2}(x))$-Laplace operator, the corresponding
integral functional, and related boundary imbedding theorem in Section 2; 
In Section 3, we state our assumptions precisely and show the existence of weak solutions to problem $(P)$ by variational arguments.

\section*{1. The spaces $W^{1,p(x)}(\Omega)$}
In this section, we will give out some theories on spaces
$W^{1,p(x)}(\Omega)$ which we call generalized Lebesgue-Sobolev
spaces. Firstly we state some basic properties of spaces
$W^{1,p(x)}(\Omega)$ which will be used later (for details see
\cite{17,14,15,4} and the references therein).

We write
\begin{align*}
&C_{+}(\overline{\Omega})=\{h|h\in C(\overline{\Omega}),h(x)>1\text{ for any }x\in\overline{\Omega}\},\\
&h^{+}=\mathop\text{max}\limits_{\overline{\Omega}}h(x),h^{-}=\mathop\text{min}\limits_{\overline{\Omega}}h(x)\text{ for any }h\in C({\overline{\Omega}}),\\
&L^{p(x)}(\Omega)=\{u|u\text{ is a measurable real-valued function, }\int_{\Omega}|u(x)|^{p(x)}\text{d}x<\infty\}.
\end{align*}
The linear vector space $L^{p(x)}(\Omega)$ can be equipped by the
following norm
\begin{align*}
|u|_{p(x)}=|u|_{L^{p(x)}(\Omega)}:=\text{inf}\bigg\{\lambda>0\bigg|\int_{\Omega}\bigg|\frac{u(x)}{\lambda}\bigg|^{p(x)}\text{d}x\leq1\bigg\},
\end{align*}
then $(L^{p(x)},|\cdot|_{p(x)})$ becomes a Banach space and we call it variable exponent Lebesgue space.\\

In the following we shall collect some basic propositions concerning the variable exponent Lebesgue spaces. These propositions will be used throughout our analysis.\\\\
\textbf{Proposition 1.1} (see Fan and Zhao \cite{17} and Zhao et al. \cite{14}).\emph{ (i) The space $(L^{p(x)},|\cdot|_{p(x)})$ is a separable, uniform convex Banach space, and its conjugate space is $L^{q(x)}(\Omega)$, where $\frac{1}{p(x)}+\frac{1}{q(x)}=1$. For any $u\in L^{p(x)}$ and $v\in L^{q(x)}$, we have
\begin{align*}
\bigg|\int_{\Omega}uv\text{\emph{d}}x\bigg|\leq(\frac{1}{p^{-}}+\frac{1}{q^{-}})|u|_{p(x)}|v|_{q(x)}.
\end{align*}
(ii) If $p_{1}, p_{2}\in C_{+}(\Omega), p_{1}(x)\leq p_{2}(x)$ for any $x\in\overline{\Omega}$, then $L^{p_{2}(x)}(\Omega)\hookrightarrow L^{p_{1}(x)}(\Omega)$, and the imbedding is continuous.
}\\\\
\textbf{Proposition 1.2} (see Fan and Zhao \cite{17} and Zhao et al. \cite{15}).\emph{ If $f:\Omega\times\mathbb{R}\rightarrow\mathbb{R}$ is a Caratheodory function and satisfies
\begin{align*}
|f(x,s)|\leq a(x)+b|s|^{\frac{p_{1}(x)}{p_{2}(x)}}\text{ for any }x\in\Omega,s\in\mathbb{R},
\end{align*}
where $p_{1},p_{2}\in C_{+}(\overline{\Omega}), a(x)\in L^{p_{2}(x)}(\Omega), a(x)\geq0$ and $b\geq0$ is a constant, then the Nemytsky operator from $L^{p_{1}(x)}(\Omega)$ to $L^{p_{2}(x)}(\Omega)$ defined by $(N_{f}(u))(x)=f(x,u(x))$ is a continuous and bounded operator.
}\\\\
\textbf{Proposition 1.3} (see Fan and Zhao \cite{17} and Zhao et al. \cite{14}).\emph{ If we denote
\begin{align*}
\rho(u)=\int_{\Omega}|u|^{p(x)}\text{\emph{d}}x, \forall u\in L^{p(x)}(\Omega),
\end{align*}
then\\
(i) $|u|_{p(x)}<1(=1;>1)\Leftrightarrow\rho(u)<1(=1;>1)$;\\
(ii) $|u|_{p(x)}>1\Rightarrow |u|_{p(x)}^{p^{-}}\leq\rho(u)\leq|u|_{p(x)}^{p^{+}};\,\,|u|_{p(x)}<1\Rightarrow |u|_{p(x)}^{p^{-}}\geq\rho(u)\geq|u|_{p(x)}^{p^{+}}$;\\
(iii)
$|u|_{p(x)}\rightarrow0\Leftrightarrow\rho(u)\rightarrow0;\,\,|u|_{p(x)}\rightarrow\infty\Leftrightarrow\rho(u)\rightarrow\infty.$
}\\\\
\textbf{Proposition 1.4} (see Fan and Zhao \cite{17} and Zhao et al. \cite{14}).\emph{ If $u, u_{n}\in L^{p(x)}(\Omega), n=1,2,...,$ then the following statements are equivalent to each other:\\
(1) $\text{\emph{lim}}_{k\rightarrow\infty}|u_{k}-u|_{p(x)}=0$;\\
(2) $\text{\emph{lim}}_{k\rightarrow\infty}\rho(u_{k}-u)=0$;\\
(3) $u_{k}\rightarrow u$ in measure in $\Omega$ and $\text{\emph{lim}}_{k\rightarrow\infty}\rho(u_{k})=\rho(u)$.
}\\\\
The variable exponent Sobolev space $W^{1,p(x)}(\Omega)$ is defined
by
\begin{align*}
W^{1,p(x)}(\Omega)=\{u\in L^{p(x)}(\Omega)| | \nabla u|\in L^{p(x)}(\Omega)\}
\end{align*}
and it is equipped with the norm
\begin{align*}
\|u\|_{p(x)}=|u|_{p(x)}+|\nabla u|_{p(x)},\forall u\in W^{1,p(x)}(\Omega).
\end{align*}
If we denote
\begin{align*}
p^{*}(x)=\left\{ \begin{array}{rcl}
         &\frac{Np(x)}{N-p(x)}, &  p(x)<N; \\
         &\infty,&p(x)\geq N.
          \end{array}\right.
\end{align*}
We have the following\\\\
\textbf{Proposition 1.5} (see Fan and Zhao \cite{17}).\emph{
(i) $W^{1,p(x)}(\Omega)$ is a separable reflexive Banach space;\\
(ii) If $q\in C_{+}(\overline{\Omega})$ and $q(x)<p^{*}(x)$ for any $x\in \overline{\Omega}$, then the embedding from $W^{1,p(x)}(\Omega)$ to $L^{q(x)}(\Omega)$ is compact and continuous.
}\\\\
\textbf{Proposition 1.6} (see Yao \cite{13}). \emph{
If we denote
\begin{align*}
p_{*}(x)=\left\{ \begin{array}{rcl}
         &\frac{(N-1)p(x)}{N-p(x)}, &  p(x)<N; \\
         &\infty,&p(x)\geq N,
          \end{array}\right.
\end{align*}
then the embedding from $W^{1,p(x)}(\Omega)$ to $L^{q(x)}(\partial\Omega)$ is compact and continuous, where $q(x)\in C_{+}(\partial\Omega)$ and $q(x)<p_{*}(x)$, $\forall x\in\partial\Omega$.
}\\

\section*{2. Properties of $(p_{1}(x),p_{2}(x))$-Laplace operator}

In this section we give the properties of the $(p_{1}(x),p_{2}(x))$-Laplace operator $(-\Delta_{p_{1}(x)}-\Delta_{p_{2}(x)})u:=-\text{div}(|\nabla u|^{p_{1}(x)-2}\nabla u)-\text{div}(|\nabla u|^{p_{2}(x)-2}\nabla u)$. Consider the following functional,
\begin{align*}
J(u)=\int_{\Omega}\frac{1}{p_{1}(x)}(|\nabla u|^{p_{1}(x)}+|u|^{p_{1}(x)})\text{d}x+\int_{\Omega}\frac{1}{p_{2}(x)}(|\nabla u|^{p_{2}(x)}+|u|^{p_{2}(x)})\text{d}x, \forall u\in X,
\end{align*}
where $X:=W^{1,p_{1}(x)}(\Omega)\cap W^{1,p_{2}(x)}(\Omega)$ with its norm given by $\|u\|:=\|u\|_{p_{1}(x)}+\|u\|_{p_{2}(x)},\forall u\in X$.

It is obvious that $J\in C^{1}(X,\mathbb{R})$ (see \cite{18}), and the $(p_{1}(x),p_{2}(x))$-Laplace operator is the derivative
operator of $J$ in the weak sense. Denote $L=J':X\rightarrow X^{*}$,
\begin{align*}
\langle L(u),v\rangle=\int_{\Omega}(|\nabla u|^{p_{1}(x)-2}&\nabla u\nabla v+|u|^{p_{1}(x)-2}uv)\text{d}x\\&+\int_{\Omega}(|\nabla u|^{p_{2}(x)-2}\nabla u\nabla v+|u|^{p_{2}(x)-2}uv)\text{d}x,\forall u,v\in X,
\end{align*}
in which $\langle\cdot,\cdot\rangle$ is the dual pair between $X$ and its dual $X^{*}$.\\\\
\textbf{Remark 2.1.} \emph{
$(X,\| \cdot \|)$ is a separable reflexive Banach space.
}\\

Similar to \cite{16}, we can deduce that\\\\
\textbf{Theorem 2.2.} \emph{
(i) $L:X\rightarrow X^{*}$ is a continuous, bounded and strictly monotone operator;\\
(ii) $L$ is a mapping of type $(S_{+})$, namely: if $u_{n}\rightharpoonup u$ in $X$ and $\overline{\text{lim}}_{n\rightarrow\infty}\langle L(u_{n})-L(u),u_{n}-u\rangle\leq 0$, then $u_{n}\rightarrow u$ in $X$;\\
(iii) $L:X\rightarrow X^{*}$ is a homeomorphism.
}\\\\

\section*{3. Solutions to the equation}
In this section we will give the existence results of weak solutions
to problem $(P)$. We denote
\begin{align*}
p_{M}(x)=\text{max}\{p_{1}(x),p_{2}(x)\},\,p_{m}(x)=\text{min}\{p_{1}(x),p_{2}(x)\},\forall
x\in \Omega.
\end{align*}
It is easy to see that $p_{M}(x),p_{m}(x)\in
C_{+}(\overline{\Omega})$.

\noindent\textbf{Lemma 3.1.}\emph{ 
(i) For $q(x)\in C_+(\overline\Omega)$ such
that $q(x)<p_{M}^{*}(x)$ for any $x\in\overline\Omega$, we have
$X:=W^{1,p_{1}(x)}(\Omega)\cap
W^{1,p_{2}(x)}(\Omega)=W^{1,p_{M}(x)}(\Omega)\hookrightarrow
L^{q(x)}(\Omega)$, and the embedding is continuous and compact;\\
(ii) For $r(x)\in C_+(\partial\Omega)$ such
that $r(x)<p_{M*}(x)$ for any $x\in\partial\Omega$, we have
$X:=W^{1,p_{1}(x)}(\Omega)\cap
W^{1,p_{2}(x)}(\Omega)=W^{1,p_{M}(x)}(\Omega)\hookrightarrow
L^{r(x)}(\partial\Omega)$, and the embedding is continuous and compact.
}\\\\
\textbf{Proof.} From $X:=W^{1,p_{1}(x)}(\Omega)\cap
W^{1,p_{2}(x)}(\Omega)=W^{1,p_{M}(x)}(\Omega)\cap W^{1,p_{m}(x)}(\Omega)=W^{1,p_{M}(x)}(\Omega)$ we can get the conclusion by Proposition 1.5 and Proposition 1.6. $\square$\\

To proceed, we give the definition of weak solution to Problem $(P)$:\\

\noindent\textbf{Definition 3.2.}\emph{ We say that $u\in X$ is a
weak solution of $(P)$ if the following equality
\begin{align*}
&\int_{\Omega}(|\nabla u|^{p_{1}(x)-2}\nabla u\nabla v+|u|^{p_{1}(x)-2}uv)\text{d}x\\
&+\int_{\Omega}(|\nabla u|^{p_{2}(x)-2}\nabla u\nabla v+|u|^{p_{2}(x)-2}uv)\text{d}x=\\
&\int_{\Omega}\lambda f(x,u)v\text{d}x+\int_{\partial\Omega}\mu g(x,u)v\text{d}\sigma
\end{align*}
holds for for any $v\in X:=W^{1,p_{1}(x)}(\Omega)\cap
W^{1,p_{2}(x)}(\Omega)$, where \emph{d}$\sigma$ is the surface measure on $\partial\Omega$.
}\\

We denote the functional $\varphi: X\rightarrow \mathbb{R}$ by
\begin{align*}
\varphi(u):=\int_{\Omega}\frac{1}{p_{1}(x)}(|\nabla
u|^{p_{1}(x)}+|u|^{p_{1}(x)})\text{d}x+\int_{\Omega}\frac{1}{p_{2}(x)}(|\nabla
u|^{p_{2}(x)}+|u|^{p_{2}(x)})\text{d}x\\
-\int_{\Omega}\lambda F(x,u)\text{d}x-\int_{\partial\Omega}\mu G(x,u)\text{d}\sigma, \forall u\in X,
\end{align*} 
where
$F(x,t)=\int^{t}_{0}f(x,s)ds$, and $G(x,t)=\int^{t}_{0}g(x,s)ds$. Next, we state the assumptions on $f$ and $g$:\\
($\text{f}_{0}$) $f:\Omega\times\mathbb{R}\rightarrow\mathbb{R}$ satisfies the Caratheodory condition and there exist two constants $C_{1}\geq0,C_{2}>0$ such that:
\begin{align*}
|f(x,t)|\leq C_{1}+C_{2}|t|^{\alpha(x)-1}, \forall (x,t)\in\Omega\times\mathbb{R},
\end{align*} 
where $\alpha(x)\in C_{+}(\overline{\Omega})$ and $\alpha(x)<p_{M}^{*}(x),\forall x\in \overline{\Omega}$. \\
($\text{f}_{1}$) There exists constants $M_{1}>0, \theta_{1}>p_{M}^{+}$ such that:
\begin{align*}
0<\theta_{1}F(x,t)\leq f(x,t)t, |t|\geq M_{1},\forall x\in\Omega.
\end{align*}
($\text{f}_{2}$) $f(x,t)=o(|t|^{p_{M}^{+}-1})$, $t\rightarrow 0$ for $x\in\Omega$ uniformly.\\
($\text{f}_{3}$) $f(x,-t)=-f(x,t)$, $\forall x\in\Omega,t\in\mathbb{R}$.\\
($\text{g}_{0}$) $g:\partial\Omega\times\mathbb{R}\rightarrow\mathbb{R}$ satisfies the Caratheodory condition and there exist two constants $C_{1}'\geq0,C_{2}'>0$ such that:
\begin{align*}
|g(x,t)|\leq C_{1}'+C_{2}'|t|^{\beta(x)-1}, \forall (x,t)\in\partial\Omega\times\mathbb{R},
\end{align*} 
where $\beta(x)\in C_{+}(\partial\Omega)$ and $\beta(x)<p_{M*}(x),\forall x\in \partial\Omega$.\\
($\text{g}_{1}$) There exists a constant $M_{2}>0, \theta_{2}>p_{M}^{+}$ such that:
\begin{align*}
0<\theta_{2}G(x,t)\leq g(x,t)t, |t|\geq M_{2},\forall x\in\partial\Omega.
\end{align*}
($\text{g}_{2}$) $g(x,t)=o(|t|^{p_{M}^{+}-1})$, $t\rightarrow 0$ for $x\in\partial\Omega$ uniformly.\\
($\text{g}_{3}$) $g(x,-t)=-g(x,t)$, $\forall x\in\partial\Omega,t\in\mathbb{R}$.

Under the assumptions ($\text{f}_{0}$) and ($\text{g}_{0}$) the functional $\varphi$ defined before is of class $C^{1}(X,\mathbb{R})$, and 
\begin{align*}
&\langle \varphi(u),v\rangle=\int_{\Omega}(|\nabla u|^{p_{1}(x)-2}\nabla u\nabla v+|u|^{p_{1}(x)-2}uv)\text{d}x\\
&+\int_{\Omega}(|\nabla u|^{p_{2}(x)-2}\nabla u\nabla v+|u|^{p_{2}(x)-2}uv)\text{d}x,\\
&-\int_{\Omega}\lambda f(x,u)v\text{d}x-\int_{\partial\Omega}\mu g(x,u)v\text{d}\sigma,\forall u,v\in X.
\end{align*}
Therefore, the weak solution of $(P)$ is
exactly the critical point of $\varphi$. Moreover, $\varphi$ is even if $f$ and $g$ are both odd with respect to the second argument respectively.\\\\
\textbf{Theorem 3.3.}\emph{ If ($\text{\emph{f}}_{0}$) and ($\text{\emph{g}}_{0}$) hold and $\alpha^{+},\beta^{+}<p_{m}^{-}$, then $(P)$ has a weak solution.
}\\\\
\textbf{Proof.} By Conditions ($\text{f}_{0}$) and ($\text{g}_{0}$) we get $|F(x,t)|\leq C(1+|t|^{\alpha(x)}),(x,t)\in\Omega\times\mathbb{R}$ and $|G(x,t)|\leq C(1+|t|^{\beta(x)}),(x,t)\in\partial\Omega\times\mathbb{R}$. And for $\|u\|$ big enough, we have
\begin{align*}
\varphi(u)&=\int_{\Omega}\frac{1}{p_{1}(x)}(|\nabla
u|^{p_{1}(x)}+|u|^{p_{1}(x)})\text{d}x+\int_{\Omega}\frac{1}{p_{2}(x)}(|\nabla
u|^{p_{2}(x)}+|u|^{p_{2}(x)})\text{d}x-\\
&\int_{\Omega}\lambda F(x,u)\text{d}x-\int_{\partial\Omega}\mu G(x,u)\text{d}\sigma\\
&\geq \frac{1}{p_{M}^{+}}(\|u_{n}\|_{p_{1}(x)}^{p_{1}^{-}}+\|u_{n}\|_{p_{2}(x)}^{p_{2}^{-}})-C|\lambda|\int_{\Omega}|u|^{\alpha(x)}\text{d}x-C|\mu|\int_{\partial\Omega}|u|^{\beta(x)}\text{d}\sigma-C_{3}\\
&\geq \frac{C_{4}}{p_{M}^{+}}\|u\|^{p_{m}^{-}}-C|\lambda|\|u\|^{\alpha^{+}}-C|\mu|\|u\|^{\beta^{+}}-C_{3}\rightarrow\infty,\text{ as }\|u\|\rightarrow\infty.
\end{align*}
So $\varphi$ is coercive. In view that $\varphi$ is also weakly lower semicontinuous, we see that $\varphi$ has a global minimum point $u\in X$, which is a weak solution to Problem $(P)$. We now complete the proof. $\square$\\\\
\textbf{Definition 3.4.} \emph{We say that the function $\varphi\in C^{1}(X,\mathbb{R})$ satisfies the Palais-Smale $(PS)$ condition in $X$ if any sequence $\{u_{n}\}\in X$ such that
\begin{align*}
|\varphi(u_{n})|\leq B,\text{ for some }B\in\mathbb{R};\\
\varphi'(u_{n})\rightarrow 0\text{ in
}X^{*}\text{ as }n\rightarrow\infty
\end{align*}
has a convergent subsequence.}\\\\
\textbf{Lemma 3.5.} \emph{
If ($\text{\emph{f}}_{0}$),($\text{\emph{f}}_{1}$),($\text{\emph{g}}_{0}$),($\text{\emph{g}}_{1}$) hold and $\lambda,\mu\geq0$, then $\varphi$ satisfies (PS) condition.}\\\\
\textbf{Proof.} Suppose that $\{u_{n}\}\subset X$, $\{\varphi(u_{n})\}$ is bounded and $\|\varphi'(u_{n})\|\rightarrow0$. Denote $\Phi(u)=-\lambda\int_{\Omega}F(x,u)\text{d}x$ and $\Psi(u)=-\mu\int_{\partial\Omega}G(x,u)\text{d}\sigma$, they are both weakly continuous and their derivative operators are compact. By Lemma 3.1, we deduce that $\varphi'=L+\Phi'+\Psi'$ is also of type ($S_{+}$). We only need to verify that $\{u_{n}\}$ is bounded. For $\|u\|$ big enough, we have
\begin{align*}
C&\geq\varphi(u_{n})=\int_{\Omega}\frac{1}{p_{1}(x)}(|\nabla
u|^{p_{1}(x)}+|u|^{p_{1}(x)})\text{d}x+\int_{\Omega}\frac{1}{p_{2}(x)}(|\nabla
u|^{p_{2}(x)}+|u|^{p_{2}(x)})\text{d}x\\
&-\int_{\Omega}\lambda F(x,u)\text{d}x-\int_{\partial\Omega}\mu G(x,u)\text{d}\sigma\\
&\geq \int_{\Omega}\frac{1}{p_{1}(x)}(|\nabla u|^{p_{1}(x)}+|u|^{p_{1}(x)})\text{d}x+\int_{\Omega}\frac{1}{p_{2}(x)}(|\nabla u|^{p_{2}(x)}+|u|^{p_{2}(x)})\text{d}x\\
&-\int_{\Omega}\frac{u_{n}}{\theta}\lambda f(x,u_{n})\text{d}x-\int_{\partial\Omega}\frac{u_{n}}{\theta}\mu g(x,u_{n})\text{d}\sigma-c\\
&\geq \int_{\Omega}(\frac{1}{p_{1}(x)}-\frac{1}{\theta})(|\nabla
u|^{p_{1}(x)}+|u|^{p_{1}(x)})\text{d}x+\int_{\Omega}(\frac{1}{p_{2}(x)}-\frac{1}{\theta})(|\nabla
u|^{p_{2}(x)}+|u|^{p_{2}(x)})\text{d}x\\
&+\frac{1}{\theta}\big(\int_{\Omega}(|\nabla u_{n}|^{p_{1}(x)}+|u_{n}|^{p_{1}(x)}+|\nabla u_{n}|^{p_{2}(x)}+|u_{n}|^{p_{2}(x)})\text{d}x-\int_{\Omega}\lambda f(x,u_{n})u_{n}\text{d}x-\int_{\partial\Omega}\mu g(x,u_{n})u_{n}\text{d}\sigma\big)-c\\
&\geq (\frac{1}{p_{M}^{+}}-\frac{1}{\theta})(\|u_{n}\|_{p_{1}(x)}^{p_{1}^{-}}+\|u_{n}\|_{p_{2}(x)}^{p_{2}^{-}})+\frac{1}{\theta}\big(\int_{\Omega}(|\nabla u_{n}|^{p_{1}(x)}+|u_{n}|^{p_{1}(x)}+|\nabla u_{n}|^{p_{2}(x)}+|u_{n}|^{p_{2}(x)})\text{d}x\\
&-\int_{\Omega}\lambda f(x,u_{n})u_{n}\text{d}x-\int_{\partial\Omega}\mu g(x,u_{n})u_{n}\text{d}\sigma\big)-c\\
&\geq C_{4}(\frac{1}{p_{M}^{+}}-\frac{1}{\theta})\|u_{n}\|^{p_{m}^{-}}-\frac{1}{\theta}\|\varphi'(u_{n})\|\|u_{n}\|-c,
\end{align*}
in which $\theta=\text{min}\{\theta_{1},\theta_{2}\}$. The above inequality implies $\{u_{n}\}$ is bounded in $X$. Finally we get that $\varphi$ satisfies (PS) condition. $\square$\\\\
\textbf{Theorem 3.6.} \emph{
If ($\text{\emph{f}}_{0}$), ($\text{\emph{f}}_{1}$), ($\text{\emph{f}}_{2}$), ($\text{\emph{g}}_{0}$), ($\text{\emph{g}}_{1}$), ($\text{\emph{g}}_{2}$) hold; $\alpha_{-},\beta_{-}>p_{M}^{+}$ and $\lambda,\mu\geq0$ such that $\lambda^{2}+\mu^{2}\not=0$, then $(P)$ has a nontrivial solution.}\\\\
\textbf{Proof.} We shall show that $\varphi$ satisfies conditions of
Mountain Pass Theorem. \\(i) From Lemma 3.5, $\varphi$ satisfies
$(PS)$ condition in $X$. Since
$p^{+}<\alpha^{-}\leq\alpha(x)<p^{*}(x)$, $p^{+}<\beta^{-}\leq\beta(x)<p^{*}(x)$, these imply $X\hookrightarrow\hookrightarrow L^{p_{M}^{+}}(\Omega)$, $X\hookrightarrow\hookrightarrow L^{p_{M}^{+}}(\partial\Omega)$, i.e. there exists $C_{0}>0$ such that
\begin{align*}
|u|_{L^{p_{M}^{+}}(\Omega)}\leq C_{0}\|u\|,\forall u\in X,|u|_{L^{p_{M}^{+}}(\partial\Omega)}\leq C_{0}\|u\|,\forall u\in X.
\end{align*}
Let $\epsilon>0$ be small enough such that $\epsilon C_{0}^{p_{M}^{+}}(\lambda+\mu)\leq\frac{C_{5}}{2p_{M}^{+}}$. By the assumptions ($\text{f}_{0}$), ($\text{f}_{2}$) and ($\text{g}_{0}$), ($\text{g}_{2}$), we have
\begin{align*}
F(x,t)\leq \epsilon|t|^{p_{M}^{+}}+C(\epsilon)|t|^{\alpha(x)},\forall(x,t)\in\Omega\times\mathbb{R};\\
G(x,t)\leq \epsilon|t|^{p_{M}^{+}}+C(\epsilon)|t|^{\beta(x)},\forall(x,t)\in\partial\Omega\times\mathbb{R}.
\end{align*}
For $\|u\|\leq 1$ we have the following
\begin{align*}
\varphi(u)&=\int_{\Omega}\frac{1}{p_{1}(x)}(|\nabla
u|^{p_{1}(x)}+|u|^{p_{1}(x)})\text{d}x+\int_{\Omega}\frac{1}{p_{2}(x)}(|\nabla
u|^{p_{2}(x)}+|u|^{p_{2}(x)})\text{d}x\\
&-\int_{\Omega}\lambda F(x,u)\text{d}x-\int_{\partial\Omega}\mu G(x,u)\text{d}\sigma\\
&\geq \frac{1}{p_{M}^{+}}(\|u\|_{p_{1}(x)}^{p_{1}^{+}}+\|u\|_{p_{2}(x)}^{p_{2}^{+}})-\lambda\epsilon\int_{\Omega}|u|^{p_{M}^{+}}\text{d}x-\lambda C(\epsilon)\int_{\Omega}|u|^{\alpha(x)}\text{d}x\\
&-\mu\epsilon\int_{\partial\Omega}|u|^{p_{M}^{+}}\text{d}\sigma-\mu C(\epsilon)\int_{\partial\Omega}|u|^{\beta(x)}\text{d}\sigma\\
&\geq \frac{C_{5}}{p_{M}^{+}}\|u\|^{p_{M}^{+}}-(\lambda+\mu)\epsilon C_{0}^{p_{M}^{+}}\|u\|^{p_{M}^{+}}- C(\epsilon)(\lambda\|u\|^{\alpha^{-}}+\mu\|u\|^{\beta^{-}})\\
&\geq \frac{C_{5}}{2p_{M}^{+}}\|u\|^{p_{M}^{+}}-C(\epsilon)(\lambda\|u\|^{\alpha^{-}}+\mu\|u\|^{\beta^{-}}),
\end{align*}
which implies the existence of $r\in (0,\,1)$ and $\delta>0$ such that $\varphi(u)\geq\delta>0$ for every $u\in X$ satisfies $\|u\|=r$.\\
(ii) From $(\text{f}_{1})$ and $(\text{g}_{1})$ we see
\begin{align*}
F(x,t)\geq C|t|^{\theta_{1}},\forall x\in\overline{\Omega}, t\geq M,\\
G(x,t)\geq C|t|^{\theta_{2}},\forall x\in\partial\Omega, t\geq M.
\end{align*}
For any fixed $w\in X\backslash\{0\}$, and $t>1$, we have
\begin{align*}
\varphi(tw)&=\int_{\Omega}\frac{1}{p_{1}(x)}(|\nabla
tw|^{p_{1}(x)}+|tw|^{p_{1}(x)})\text{d}x+\int_{\Omega}\frac{1}{p_{2}(x)}(|\nabla
tw|^{p_{2}(x)}+|u|^{p_{2}(x)})\text{d}x\\
&-\int_{\Omega}\lambda F(x,tw)\text{d}x-\int_{\partial\Omega}\mu G(x,tw)\text{d}\sigma\\
&\leq t^{p_{M}^{+}}(\int_{\Omega}\frac{1}{p_{1}(x)}(|\nabla
w|^{p_{1}(x)}+|w|^{p_{1}(x)})\text{d}x+\int_{\Omega}\frac{1}{p_{2}(x)}(|\nabla
w|^{p_{2}(x)}+|w|^{p_{2}(x)})\text{d}x)\\
&-C\lambda t^{\theta_{1}}\int_{\Omega}|w|^{\theta_{1}}\text{d}x-C\mu t^{\theta_{2}}\int_{\partial\Omega}|w|^{\theta_{2}}\text{d}x-C_{6}\\
&\rightarrow-\infty\text{ as } t\rightarrow+\infty,\text{ since }\theta_{1},\theta_{2}>p_{M}^{+}.
\end{align*}
(iii) It is obvious $\varphi(0)=0$ by ($\text{f}_{2}$) and ($\text{g}_{2}$).

From (i), (ii) and (iii), we conclude $\varphi$ satisfies the conditions of Mountain Pass Theorem (see \cite{18}). So $\varphi$ admits at least one nontrivial critical point. $\square$\\\\
\textbf{Theorem 3.7.} \emph{
If ($\text{\emph{f}}_{0}$), ($\text{\emph{f}}_{1}$), ($\text{\emph{f}}_{3}$), ($\text{\emph{g}}_{0}$), ($\text{\emph{g}}_{1}$), ($\text{\emph{g}}_{3}$) hold and $\max\{\alpha^{+}, \beta^{+}\}>p_{M}^{+}$; $\lambda,\mu>0$, then $\varphi$ has a sequence of critical points $\{u_{n}\}$ such that $\varphi(u_{n})\rightarrow+\infty$ and $(P)$ has infinite many pairs of solutions.}\\\\
Because $X$ is a reflexive and separable Banach space, there are $\{e_{j}\}\subset X$ and $\{e^{*}_{j}\}\subset X^{*}$ such that
\begin{align*}
X=\overline{\text{span}\{e_{j}|j=1,2,...\}},X^{*}=\overline{\text{span}\{e^{*}_{j}|j=1,2,...\}}
\end{align*}
and
\begin{align*}
\langle e_{i}^{*}, e_{j}\rangle=\left\{ \begin{array}{rcl}
         &1,&  i=j; \\
         &0,&i\neq j.
          \end{array}\right.
\end{align*}
Denote $X_{j}=\text{span}\{e_{j}\},Y_{k}=\oplus_{j=1}^{k}X_{j},Z_{k}=\oplus_{j=k}^{\infty}X_{j}$.\\\\
\textbf{Lemma 3.8.} \emph{
If $\alpha\in C_{+}(\overline{\Omega}),\alpha(x)<p_{M}^{*}(x)$ for any $x\in\overline{\Omega}$; and $\beta\in C_{+}(\partial\Omega),\beta(x)<p_{M*}(x)$ for any $x\in\partial\Omega$, denote
\begin{align*}
\alpha_{k}=\text{sup}\{|u|_{L^{\alpha(x)}(\Omega)}|\|u\|=1,u\in Z_{k}\},\\
\beta_{k}=\text{sup}\{|u|_{L^{\beta(x)}(\partial\Omega)}|\|u\|=1,u\in Z_{k}\}.
\end{align*}
Then $\text{\emph{lim}}_{k\rightarrow\infty}\alpha_{k}=0$, $\text{\emph{lim}}_{k\rightarrow\infty}\beta_{k}=0$.}\\\\
\textbf{Proof.} For $0<\alpha_{k+1}\leq\alpha_{k}$, then $\alpha_{k}\rightarrow\alpha\geq0$. Suppose $u_{k}\in Z_{k}$ satisfy
\begin{align*}
\|u_{k}\|=1,0\leq\alpha_{k}-|u_{k}|_{\alpha(x)}<\frac{1}{k},
\end{align*}
then we may assume up a subsequence that $u_{k}\rightharpoonup u$ in
$X$, and
\begin{align*}
\langle e^{*}_{j},u\rangle=\text{lim}_{k\rightarrow\infty}\langle e^{*}_{j},u_{k}\rangle=0,j=1,2,...,
\end{align*}
in which the last equal sign holds for $u_{k}\in Z_{k}$. The above equality implies that $u=0$, so $u_{k}\rightharpoonup 0$ in $X$. Since the imbedding from $X$ to $L^{\alpha(x)}(\Omega)$ is compact, we have $u_{k}\rightarrow 0$ in $L^{\alpha(x)}(\Omega)$. We finally get $\alpha_{k}\rightarrow0$. The result for $\beta_{k}$ can be obtained  by the same procedure. $\square$\\\\
\textbf{Proof of Theorem 3.7.} By $(\text{f}_{1}), (\text{f}_{3})$, $(\text{g}_{1}),(\text{g}_{3})$ $\varphi$ is an even functional and satisfies $(PS)$ condition. We only need to prove that if $k$ is large enough, then there exist $\rho_{k}>\gamma_{k}>0$ such that\\
$(A_{1})$ $b_{k}:=\text{inf}\{\varphi(u)|u\in Z_{k},\|u\|=\gamma_{k}\}\rightarrow\infty,(k\rightarrow\infty)$;\\
$(A_{2})$ $a_{k}:=\text{max}\{\varphi(u)|u\in Y_{k},\|u\|=\rho_{k}\}\leq0$,\\
then the conclusion of Theorem 3.7 can be obtained by Fountain Theorem (see \cite{19}) and Lemma 3.5.\\
$(A_{1})$ For any $u\in Z_{k}$ with $\|u\|$ is big enough, we have
\begin{align*}
\varphi(u)&=\int_{\Omega}\frac{1}{p_{1}(x)}(|\nabla
u|^{p_{1}(x)}+|u|^{p_{1}(x)})\text{d}x+\int_{\Omega}\frac{1}{p_{2}(x)}(|\nabla
u|^{p_{2}(x)}+|u|^{p_{2}(x)})\text{d}x\\
&-\int_{\Omega}\lambda F(x,u)\text{d}x-\int_{\partial\Omega}\mu G(x,u)\text{d}\sigma\\
&\geq \frac{1}{p_{M}^{+}}(\|u\|_{p_{1}(x)}^{p_{1}^{-}}+\|u\|_{p_{2}(x)}^{p_{2}^{-}})-c\lambda\int_{\Omega}|u|^{\alpha(x)}\text{d}x-c\mu\int_{\partial\Omega}|u|^{\beta(x)}\text{d}\sigma-c_{1}\\
&\geq \frac{C_{4}}{p_{M}^{+}}\|u\|^{p_{m}^{-}}-c\lambda|u|_{L^{\alpha(x)}(\Omega)}^{\alpha(\xi)}-c\mu|u|_{L^{\beta(x)}(\partial\Omega)}^{\beta(\eta)}-c_{2},\text{ where }\xi\in\Omega,\eta\in\partial\Omega\\
&\geq\left\{ \begin{array}{rcl}
         &\frac{C_{4}}{p_{M}^{+}}\|u\|^{p_{m}^{-}}-c-c_{2},&\text{ if }|u|_{\alpha(x)}\leq 1,|u|_{\beta(x)}\leq 1; \\
         &\frac{C_{4}}{p_{M}^{+}}\|u\|^{p_{m}^{-}}-c\alpha_{k}^{\alpha^{+}}\|u\|^{\alpha^{+}}-c'_{2},&\text{ if }|u|_{\alpha(x)}> 1,|u|_{\beta(x)}\leq 1;\\
         &\frac{C_{4}}{p_{M}^{+}}\|u\|^{p_{m}^{-}}-c\beta_{k}^{\beta^{+}}\|u\|^{\beta^{+}}-c''_{2},&\text{ if }|u|_{\alpha(x)}\leq 1,|u|_{\beta(x)}> 1; \\
         &\frac{C_{4}}{p_{M}^{+}}\|u\|^{p_{m}^{-}}-c\alpha_{k}^{\alpha^{+}}\|u\|^{\alpha^{+}}-c\beta_{k}^{\beta^{+}}\|u\|^{\beta^{+}}-c_{2},&\text{ if }|u|_{\alpha(x)}> 1,|u|_{\beta(x)}> 1; \\
          \end{array}\right.\\
&\geq \frac{C_{4}}{p_{M}^{+}}\|u\|^{p_{m}^{-}}-c\alpha_{k}^{\alpha^{+}}\|u\|^{\alpha^{+}}-c\beta_{k}^{\beta^{+}}\|u\|^{\beta^{+}}-c_{3}\\
&=
C_{4}(\frac{1}{p_{M}^{+}}\|u\|^{p_{m}^{-}}-c'\alpha_{k}^{\alpha^{+}}\|u\|^{\alpha^{+}}-c'\beta_{k}^{\beta^{+}}\|u\|^{\beta^{+}})-c_{3}\\
&\geq
C_{4}(\frac{1}{p_{M}^{+}}\|u\|^{p_{m}^{-}}-c''\zeta_{k}^{\zeta^{+}}\|u\|^{\zeta^{+}})-c_{3},
\end{align*}
in which $\zeta^{+}=\text{max}\{\alpha^{+},\beta^{+}\},\zeta_{k}=\text{max}\{\alpha_{k},\beta_{k}\}$.
Set
$\|u\|=\gamma_{k}=(c''\zeta^{+}\zeta_{k}^{\zeta^{+}})^{\frac{1}{p_{m}^{-}-\zeta^{+}}}$.
Because $\zeta_{k}\rightarrow 0$ and $p_{m}^{-}\leq
p_{M}^{+}<\zeta^{+}$, we have
\begin{align*}
\varphi(u)&\geq C_{4}(\frac{1}{p_{M}^{+}}\|u\|^{p_{m}^{-}}-c''\zeta_{k}^{\zeta^{+}}\|u\|^{\zeta^{+}})-c_{3}\\
&=C_{4}(\frac{1}{p_{M}^{+}}(c''\zeta^{+}\zeta_{k}^{\zeta^{+}})^{\frac{p_{m}^{-}}{p_{m}^{-}-\zeta^{+}}}-c''\zeta_{k}^{\alpha^{+}}(c''\zeta^{+}\zeta_{k}^{\zeta^{+}})^{\frac{\zeta^{+}}{p_{m}^{-}-\zeta^{+}}})-c_{3}\\
&=C_{4}(\frac{1}{p_{M}^{+}}-\frac{1}{\zeta^{+}})(c''\zeta^{+}\zeta_{k}^{\zeta^{+}})^{\frac{p_{m}^{-}}{p_{m}^{-}-\zeta^{+}}}-c_{3}\rightarrow\infty, \text{ as }k\rightarrow\infty.
\end{align*}
$(A_{2})$ From $(\text{f}_{1})$, we get
\begin{align*}
F(x,t)\geq c_{1}|t|^{\theta_{1}}-c_{2},\forall (x,t)\in\Omega\times\mathbb{R};\\
G(x,t)\geq c_{1}|t|^{\theta_{2}}-c_{2},\forall (x,t)\in\partial\Omega\times\mathbb{R}.
\end{align*}
Because $\theta_{1}>p_{M}^{+}$ and $\text{dim}Y_{k}=k<\infty$ (all norms are equivalent in $Y_{k}$), it is easy to get $\varphi(u)\rightarrow-\infty$ as $\|u\|\rightarrow\infty$ for $u\in Y_{k}$. $\square$\\\\
\textbf{Theorem 3.9.} \emph{
Let $\alpha(x)\in C_{+}(\overline{\Omega})$, $\beta\in C_{+}(\partial\Omega)$ and $\alpha(x)<p^{*}(x),\forall x\in\overline{\Omega}$; $\beta(x)<p_{M*}(x), \forall x\in\partial\Omega$. If $f(x,t)=|t|^{\alpha(x)-2}t$, $g(x,t)=|t|^{\beta(x)-2}t$, $\alpha^{-}>p_{M}^{+}$, $\beta^{+}<p_{m}^{-}$, then we have\\
(i) For every $\lambda>0$, $\mu\in\mathbb{R}$, $(P)$ has a sequence of weak solutions $\{\pm u_{k}\}$ such that $\varphi(\pm u_{k})\rightarrow\infty, k\rightarrow \infty$;\\
(ii) For every $\mu>0$, $\lambda\in\mathbb{R}$, $(P)$ has a sequence of weak solutions $\{\pm v_{k}\}$ such that $\varphi(\pm v_{k})<0$, and $\varphi(\pm v_{k})\rightarrow0, k\rightarrow \infty$.
}\\\\
\textbf{Proof.} (i) The proof is similar to that of Theorem 3.7 if we use the Fountain Theorem, we only verify $(PS)$ condition here. We only need to verify the $(PS)$ sequence $\{u_{n}\}$ is bounded in $X$ as in Lemma 3.5. We assume
\begin{align*}
\{u_{n}\}\subset X, \varphi(u_{n})\leq M, \varphi'(u_{n})\rightarrow0,\text{ as }n\rightarrow\infty.
\end{align*}
For $\|u_{n}\|$ big enough,
\begin{align*}
C&\geq\varphi(u_{n})=[\int_{\Omega}\frac{1}{p_{1}(x)}(|\nabla
u_{n}|^{p_{1}(x)}+|u_{n}|^{p_{1}(x)})\text{d}x+\int_{\Omega}\frac{1}{p_{2}(x)}(|\nabla
u_{n}|^{p_{2}(x)}+|u_{n}|^{p_{2}(x)})\text{d}x\\
&-\lambda\int_{\Omega}\frac{1}{\alpha(x)}|u_{n}|^{\alpha(x)}\text{d}x-\mu\int_{\partial\Omega}\frac{1}{\beta(x)}|u_{n}|^{\beta(x)}\text{d}\sigma]\\
&-\frac{1}{\alpha^{-}}[\int_{\Omega}(|\nabla
u_{n}|^{p_{1}(x)}+|u_{n}|^{p_{1}(x)})\text{d}x+\int_{\Omega}(|\nabla
u_{n}|^{p_{2}(x)}+|u_{n}|^{p_{2}(x)})\text{d}x\\
&-\lambda\int_{\Omega}|u_{n}|^{\alpha(x)}\text{d}x-\mu\int_{\partial\Omega}|u_{n}|^{\beta(x)}\text{d}\sigma]+\frac{1}{\alpha^{-}}\langle\varphi'(u_{n}),u_{n}\rangle\\
&\geq (\frac{1}{p_{M}^{+}}-\frac{1}{\alpha^{-}})(\|u_{n}\|_{p_{1}(x)}^{p_{1}^{-}}+\|u_{n}\|_{p_{2}(x)}^{p_{2}^{-}})-C\|u_{n}\|^{\beta^{+}}-\frac{1}{\alpha^{-}}\|\varphi'(u_{n})\|_{X^{*}}\|u_{n}\|\\
&\geq C_{4}(\frac{1}{p_{M}^{+}}-\frac{1}{\alpha^{-}})\|u_{n}\|^{p_{m}^{-}}-C\|u_{n}\|^{\beta^{+}}-\frac{1}{\alpha^{-}}\|u_{n}\|.
\end{align*}
Since $\alpha^{-}>p_{m}^{+}, \beta^{+}<p_{m}^{-}$, we have $\{u_{n}\}$ is bounded in $X$.\\
(ii) We use the Dual Fountain Theorem (see \cite{19}) to prove (ii), i.e. we have to prove that $\varphi$ satisfies the \textbf{$(PS)_{c}^{*}$ condition} and there exist $\rho_{k}>r_{k}>0$ such that for large $k$,\\
$(B_{1})$ $a_{k}:=\inf\{\varphi(u)|u\in Z_{k},\|u\|=\rho_{k}\}\geq0$;\\
$(B_{2})$ $b_{k}:=\max\{\varphi(u)|u\in Y_{k},\|u\|=r_{k}\}<0$;\\
$(B_{3})$ $d_{k}:=\inf\{\varphi(u)|u\in Z_{k},\|u\|\leq\rho_{k}\}\rightarrow0,k\rightarrow\infty$.\\
We will do this step by step:\\
\textbf{$(\text{B}_{2})$} For $u\in Y_{k}$, $\|u\|$ small enough,
\begin{align*}
\varphi(u)&=\int_{\Omega}\frac{1}{p_{1}(x)}(|\nabla
u|^{p_{1}(x)}+|u|^{p_{1}(x)})\text{d}x+\int_{\Omega}\frac{1}{p_{2}(x)}(|\nabla
u|^{p_{2}(x)}+|u|^{p_{2}(x)})\text{d}x\\
&-\lambda\int_{\Omega}\frac{1}{\alpha(x)}|u|^{\alpha(x)}\text{d}x-\mu\int_{\partial\Omega}\frac{1}{\beta(x)}|u|^{\beta(x)}\text{d}\sigma\\
&\leq \frac{C_{6}}{p_{m}^{-}}\|u\|^{p^{-}}+\frac{|\lambda|}{\alpha^{-}}\int_{\Omega}|u|^{\alpha(x)}\text{d}x-\frac{\mu}{\beta^{+}}\int_{\partial\Omega}|u|^{\beta(x)}\text{d}\sigma.
\end{align*}
Since $\dim Y_{k}<\infty$ and $p_{m}^{-}>\beta^{+},\alpha^{-}>p_{M}^{+}>\beta^{+}$, if we choose $r_{k}$ small enough, $b_{k}<0$.\\
\textbf{$(\text{B}_{1})$} For $u\in Z_{k}$, $\|u\|$ small,
\begin{align*}
\varphi(u)&\geq\frac{C_{5}}{p_{M}^{+}}\|u\|^{p_{M}^{+}}-\frac{C|\lambda|}{\alpha^{-}}\|u\|^{\alpha^{-}}-\frac{\mu}{\beta^{-}}\int_{\partial\Omega}|u|^{\beta(x)}\text{d}\sigma\\
&\geq\frac{C_{5}}{p_{M}^{+}}\|u\|^{p_{M}^{+}}-\frac{C|\lambda|}{\alpha^{-}}\|u\|^{\alpha^{-}}-\frac{\mu}{\beta^{-}}\max\{|u|^{\beta^{+}}_{L^{\beta(x)}(\partial\Omega)},|u|^{\beta^{-}}_{L^{\beta(x)}(\partial\Omega)}\}.
\end{align*}
Since $\alpha^{-}>p_{M}^{+}$, there exists $\rho_{0}$ small enough such that $\frac{C|\lambda|}{\alpha^{-}}\|u\|^{\alpha^{-}}\leq\frac{C_{5}}{2p_{M}^{+}}$ for any $u:0<\|u\|\leq\rho_{0}$. Then
\begin{align*}
\varphi(u)&\geq\frac{C_{5}}{2p_{M}^{+}}\|u\|^{p_{M}^{+}}-\frac{\mu}{\beta^{-}}\max\{|u|^{\beta^{+}}_{L^{\beta(x)}(\partial\Omega)},|u|^{\beta^{-}}_{L^{\beta(x)}(\partial\Omega)}\}\\
&\geq\frac{C_{5}}{2p_{M}^{+}}\|u\|^{p_{M}^{+}}-\frac{\mu}{\beta^{-}}\max\{\beta_{k}^{\beta^{-}}\|u\|^{\beta^{-}},\beta_{k}^{\beta^{+}}\|u\|^{\beta^{+}}\}\\
&\geq\frac{C_{5}}{2p_{M}^{+}}\|u\|^{p_{M}^{+}}-\frac{\mu}{\beta^{-}}\beta_{k}^{\beta^{-}}\|u\|^{\beta^{-}}.
\end{align*}
Choose $\rho_{k}=(\frac{2p_{M}^{+}\mu\beta_{k}^{\beta^{-}}}{C_{5}\beta^{-}})^{\frac{1}{p_{M}^{+}-\beta^{-}}}$, then
\begin{align*}
\varphi(u)\geq\frac{C_{5}}{2p_{M}^{+}}\rho_{k}^{p_{M}^{+}}-\frac{C_{5}}{2p_{M}^{+}}\rho_{k}^{p_{M}^{+}}=0.
\end{align*}
From $p_{m}^{-}>\beta^{+},\beta_{k}\rightarrow0$, we know that $\rho_{k}\rightarrow0$ as $k\rightarrow\infty$.\\
\textbf{$(\text{B}_{3})$} From the proof of $(B_{1})$, we know for $u\in Z_{k}, \|u\|\leq\rho_{k}$ small enough
\begin{align*}
\varphi(u)\geq-\frac{\mu}{\beta^{-}}\beta_{k}^{\beta^{-}}\|u\|^{\beta^{-}}.
\end{align*}
Since $\beta_{k}\rightarrow0$ and $\rho_{k}\rightarrow0$ as $k\rightarrow\infty$. And from the above proof, we can choose $\rho_{k}>r_{k}>0$.\\
\textbf{The proof of $(PS)^{*}_{c}$ condition.} Consider a sequence $\{u_{n_{j}}\}\subset X$ such that
\begin{align*}
n_{j}\rightarrow\infty,u_{n_{j}}\in Y_{n_{j}},\varphi(u_{n_{j}})\rightarrow c, \varphi|'_{Y_{n_{j}}}(u_{n_{j}})\rightarrow0.
\end{align*}
Assume $\|u\|$ is big enough, for large $n$, we have:\\
Case (1), $\lambda\geq0$. We have
\begin{align*}
C&\geq\varphi(u_{n_{j}})= \varphi(u_{n_{j}})-\frac{1}{\alpha^{-}}\langle\varphi'(u_{n_{j}}),u_{n_{j}}\rangle+\frac{1}{\alpha^{-}}\langle\varphi'(u_{n_{j}}),u_{n_{j}}\rangle\\
&\geq C_{4}(\frac{1}{p_{M}^{+}}-\frac{1}{\alpha^{-}})\|u_{n_{j}}\|^{p_{m}^{-}}-C\|u_{n_{j}}\|^{\beta^{+}}-\frac{1}{\alpha^{-}}\|u_{n_{j}}\|.
\end{align*}
Since $\alpha^{-}>p^{+}\geq p^{-}>\beta^{+}$, we deduce $\{u_{n_{j}}\}$ is bounded in $X$.\\
Case (2), $\lambda<0$. We consider the following inequality to get that $\{u_{n_{j}}\}$ is bounded in $X$:
\begin{align*}
C\geq\varphi(u_{n_{j}})= \varphi(u_{n_{j}})-\frac{1}{\alpha^{+}}\langle\varphi'(u_{n_{j}}),u_{n_{j}}\rangle+\frac{1}{\alpha^{+}}\langle\varphi'(u_{n_{j}}),u_{n_{j}}\rangle.
\end{align*}
Going if necessary to a subsequence, we can assume that $u_{n_{j}}\rightharpoonup u$ in $X$. As $X=\overline{\cup_{n_{j}}Y_{n_{j}}}$, we can choose $v_{n_{j}}\in Y_{n_{j}}$ such that $v_{n_{j}}\rightharpoonup u$. Hence
\begin{align*}
\lim_{n_{j}\rightarrow\infty}\varphi'(u_{n_{j}})(u_{n{j}}-u)&=\lim_{n_{j}\rightarrow\infty}\varphi'(u_{n_{j}})(u_{n{j}}-v_{n_{j}})+\lim_{n_{j}\rightarrow\infty}\varphi'(u_{n_{j}})(v_{n{j}}-u)\\
&=\lim_{n_{j}\rightarrow\infty}(\varphi|_{Y_{n_{j}}})'(u_{n_{j}})(u_{n{j}}-v_{n_{j}})=0.
\end{align*}
As $\varphi'$ is of type $(S_{+})$, we can conclude $u_{n_{j}}\rightarrow u$. Furthermore, we have $\varphi'(u_{n_{j}})\rightarrow\varphi'(u)$. It remains to prove $\varphi'(u)=0$. Take any $w_{k}\in Y_{k}$, for $n_{j}\geq k$
\begin{align*}
\varphi'(u)w_{k}&=(\varphi'(u)-\varphi'(u_{n_{j}}))w_{k}+\varphi'(u_{n_{j}})w_{k}\\
&=(\varphi'(u)-\varphi'(u_{n_{j}}))w_{k}+(\varphi|_{Y_{n_{j}}})'(u_{n_{j}})w_{k}.
\end{align*}
Let $n_{j}\rightarrow\infty$ on the right hand side of the above equation, we get $\varphi'(u)w_{k}=0, \forall w_{k}\in Y_{k}$, which implies $\varphi'(u)=0$. This shows $\varphi$ satisfies the $(PS)^{*}_{c}$ condition for every $c\in \mathbb{R}$. $\square$

\textbf{Remark 3.10.} \emph{ We could extend all our results directly to the multi $p(x)-$Laplacian case to consider the following $(p_{1}(x), p_{2}(x),...,p_{n}(x))$-Laplace equation
\begin{align*}
(P_{n})\left\{ \begin{array}{rcl}
         &-\sum_{j=1}^{n}\Delta_{p_j(x)}u+a_j(x)|u|^{p_j (x)-2}u=f(x,u), &  \mbox{ in }
         \Omega; \\
         &|\nabla u|^{p_{1}(x)-2}\frac{\partial u}{\partial\nu}+...+|\nabla u|^{p_{n}(x)-2}\frac{\partial u}{\partial\nu}=\mu g(x,u),&\mbox{ on } \partial\Omega.
          \end{array}\right.
\end{align*}
}
where $a_j(x)\geq \delta >0$ are smooth functions  and $\delta$ a positive real numbers. With slight elaboration and necessary preparation, we hope to deal with cases where some $a_{j}(x)$ touches zero with positive zero level set.

\renewcommand{\baselinestretch}{0.1}

\end{document}